\documentclass[10pt]{amsart}

\usepackage{amssymb}
\usepackage[T1]{fontenc} 
\usepackage{url}
\usepackage{color}

\theoremstyle{plain}
\newtheorem{theorem}{Theorem}[section]
\newtheorem{proposition}[theorem]{Proposition}
\newtheorem{lemma}[theorem]{Lemma}
\newtheorem{corollary}[theorem]{Corollary}
\newtheorem{fact}[theorem]{Fact}

\newtheorem*{claim}{Claim}
\theoremstyle{definition}

\newtheorem{example}[theorem]{Example}

\newtheorem*{scvsascproblem*}{SC vs ASC problem}

\renewcommand{\=}{\approx}
\renewcommand{\leq}{\leqslant}

\newcommand{\meet}{\wedge}
\newcommand{\bigmeet}{\bigwedge}

\newcommand{\Con}{\operatorname{Con}}
\newcommand{\id}{\operatorname{id}}

\newcommand{\disc}{\operatorname{disc}}
\newcommand{\diagel}{{ \sf diag^{el}}}
\newcommand{\diagplus}{\sf{diag^+}}

\newcommand{\V}{\mathcal V}
\newcommand{\W}{\mathcal W}

\newcommand{\Q}{\mathcal Q}

\renewcommand{\S}{\mathbf S}
\renewcommand{\P}{\mathbf P}

\newcommand{\A}{\mathbf{A}}
\newcommand{\B}{\mathbf{B}}

\newcommand{\Fr}{\mathbf{F}}
\newcommand{\G}{\mathbf{G}}

\newcommand{\Nat}{\mathbb N}

\begin{document}

\title[On SC vs ASC problem]{On structural completeness vs almost structural completeness problem:\\ A discriminator varieties case study}

\author{Miguel Campercholi}
\address{Facultad de Matem\'atica, Astronom\'ia y F\'isica, Ciudad Universitaria, C\'ordoba, Argentina}
\email{camper@famaf.unc.edu.ar}

\author{Micha\l~M. Stronkowski}

\address{Faculty of Mathematics and Information Sciences,
Warsaw University of Technology Poland}
\email{m.stronkowski@mini.pw.edu.pl}

\author{Diego Vaggione}

\address{Facultad de Matem\'atica, Astronom\'ia y F\'isica, Ciudad Universitaria, C\'ordoba, Argentina}
\email{vaggione@famaf.unc.edu.ar}

\thanks{The work of the first and third authors was supported by Conicet and Secyt-UNC. The work of the second author was supported by the Polish National Science Centre [DEC-2011/01/D/ST1/06136]. }



\begin{abstract}
We study the following problem: Determine which almost structurally complete quasivarieties are structurally complete. We propose a general solution to this problem and then a solution in the semisimple case. As a consequence, we obtain a characterization of structurally complete discriminator varieties. 

An interesting corollary in logic follows: Let $L$ be a consistent propositional logic/deductive system in the language with formulas for \emph{verum}, which is a theorem, and \emph{falsum}, which is not a theorem. Assume also that $L$ has an adequate semantics given by a discriminator variety. Then $L$ is structurally complete if and only if it is  maximal. 
All such logics/deductive systems are almost structurally complete.
\end{abstract}

\maketitle

\section{Introduction}
A quasivariety $\Q$ is structurally complete (SC for short) if all its admissible (true in free algebras) quasi-identities are true in $\Q$ \cite{Ber88}. And $\Q$ is almost structurally complete (ASC for short) if every its admissible quasi-identities which do not hold in $\Q$ is passive \cite{DS14,MR13}. These are quasi-identities whose premises are not satisfiable in free algebras. Both notions have their origins in the theory of deductive systems. The SC property was introduced by Pogorzelski in \cite{Pog71} and the ASC property by Dzik in \cite{Dzi06}. From this perspective both of them serve the same purpose: to separate deductive systems whose ``proof-power'' cannot be strenghtened by adding new rules while keeping the set of theorems unchanged. Indeed, such strenghtening cannot be obtained by adding passive rules. A detailed discussion concerning the logical context may be found in \cite{DS14}.

It appeared that there are important deductive systems, like S5 modal logic augmented by
{\it modus ponens rule} and {\it generalization rule}, or \L{}ukasiewicz $n$-valued logic augmented by {\it modus ponens rule}, which are not SC merely because of the underivability of some passive rules. Thus they are ASC. (Note that both these examples have adequate semantics given by discriminator varieties.) One then can ask: are such examples rare, or maybe there are plenty of them (in whatever sense)? This leads us to the following task.

\begin{scvsascproblem*}
Determine which ASC qausivarieties/deductive systems are SC.
\end{scvsascproblem*}

There are two, quite obvious, conditions  for an ASC quasivarieties that yield SC: minimality and that
there is an idempotent element in a free algebra. A relaxation of the second
condition would be that every nontrivial algebra from the quasivariety under consideration admits a homomorphism into a free algebra. And indeed, it was shown in \cite[Corollary 3.5.]{DS14} that it yields SC for ASC quasivarieties. Still, there are SC quasivarieties in which this condition fails, see Examples \ref{exm:: idemp w G} and \ref{exm:: minimal discriminator without homo into F}. Nonetheless, our first result in this paper says that it is quite close to the right condition. We prove that an ASC quasivariety $\Q$ is SC iff every nontrivial algebra in $\Q$ admits a homomorphism into some elementary extension of a free algebra for $\Q$ iff every nontrivial $\Q$-finitely presented algebra admits a homomorphism into a free algebra for $\Q$ (Theorem \ref{thm::general solution}).

One could argue that the obtained condition is not applicable. On the other hand, it is hard to expect anything better on this level of generality. This is why in the subsequent considerations we restrict to semisimple quasivarieties. Here we can replace a condition with homomorphisms for the disjunction of the minimality and the existence of an idempotent element in an elementary extension of a free algebra (Theorem \ref{thm:: main}). This, in particular, yields a characterization of SC discriminator varieties (Theorem \ref{thm:: disr var}). Indeed, discriminator varieties are semisimple \cite[Theorem IV.9.4]{BS81}. Moreover, Burris proved in \cite[Theorem 3.1]{Bur92} that every discriminator variety has projective unification (see also \cite{Wro95}). This property is known to imply ASC \cite[Corollary 6]{Dzi11}\cite[Corollary 5.1]{DS14}. Discriminator varieties, introduced by Pixley in \cite{Pix70}, form a restrictive class from the perspective of universal algebra. On the other hand, algebras in discriminator varieties have a well understood structure as Boolean products of simple and trivial algebras \cite{BW77,BS81,Wer78}. (We used it to show that minimal discriminator varieties are also minimal as quasivarieties, see Proposition \ref{prop::minimal discr}). Moreover, discriminator varieties appear quite often in algebraic logic, see e.g. \cite{GFLP09,Hal62,HMT85,Jip93,KK06,Kow98,Kow04,Mad06}. In most of these cases there are distinct constants for \emph{verum} and \emph{falsum}. As a consequence, being SC is equivalent to be minimal (Corollary \ref{cor:: disc with constants}).

It is worth adding that almost all (in a certain strict sense) finite algebras in a finite language with at least one at least binary operation generate minimal discriminator varieties. They are also minimal as quasivarieties and hence are SC. It was proved by Murski{\u\i} in \cite{Mur75} (see also \cite[Chapter 6]{Ber12}) that almost all finite algebras with at least one at least binary operation have all their idempotent operations as term operations (they are idemprimal). In particular, they have the  discriminator operation as a term operation and do not have proper subalgebras.

Finally note that the SC vs ASC problem was solved in \cite[Proposition 8.6]{DS14} for varieties of closure algebras (equivalently for normal extensions of S4 modal logic) and for varieties of bounded lattices \cite[Proposition 7.7]{DS14}. Both of these cases are quite different than the case of discriminator varieties. 
It is shown there that an ASC variety of closure algebras is SC iff it does not contain a four-element simple algebra iff it satisfies the McKinsey identity. Note that a nontrivial closure algebra satisfies the McKinsey identity iff it admits a homomorphism onto a two-element algebra (a two-element closure algebra is free in every nontrivial variety of closure algebras). And the variety of bounded lattices is SC iff it is ASC iff it satisfies distributivity law.

\section{Preliminaries}
\label{sec::quasivariety theory}
In the paper we assume that the reader has basic knowledge in universal algebra \cite{BS81,Ber12} and model theory \cite{CK90,Hod93}. However, 
we recall needed concepts from quasivariety theory \cite{Gor98,Mal70} since it is less known.

A \emph{quasi-identity} is a sentence in first order algebraic language of the form
\[
(\forall \bar x)\, [s_1(\bar x)\=t_1(\bar x)\meet\cdots\meet s_{n}(\bar x)\=t_{n}(\bar x)\to s(\bar x)\=t(\bar x)],
\]
where $n\in\mathbb N$. We allow $n$ to be zero, and in such case we call the sentence an \emph{identity}. It will be convenient to have a more compact notation for quasi-identities, and we will often write them in the form
\[
(\forall \bar x)\, [\varphi(\bar x)\to \psi(\bar x)],
\]
where $\varphi$ is a conjunction of equations (i.e., atomic formulas) and $\psi$ is an equation. We call $\varphi$ the \emph{premise} and $\psi$ the \emph{conclusion} of a quasi-identity.

By the \emph{(quasi-)equational theory} of a class $\mathcal K$ of algebras in the same language we mean the set of (quasi-)identities true in $\mathcal K$.
A \emph{(quasi)variety} is a class defined by (quasi-)identities. Equivalently, a class of algebras in the same language is a quasivariety if it is closed under taking substructures, direct products and ultraproducts. If it is additionally closed under taking homomorphic images, it is a variety. (We tacitly assume that all considered classes contain algebras in the same language and are closed under taking isomorphic images. Also all considered class operators are assumed to be composed with the isomorphic image class operator.) A (quasi)variety is \emph{trivial} if it consists of one-element algebras, and is \emph{minimal} if it properly contains only the trivial (quasi)variety. We say that a (quasi)variety $\Q$ is \emph{generated by} a class $\mathcal K$ if it is the smallest (quasi)variety containing $\mathcal K$. This means that $\Q$ is defined by the \mbox{(quasi-)equational} theory of $\mathcal K$. We denote such (quasi)variety by ${\sf V}(\mathcal K)$ (${\sf Q}(\mathcal K)$ respectively). When $\mathcal K=\{\A\}$ we simplify the notation by writing ${\sf V}(\A)$ (${\sf Q}(\A)$ respectively).
Let $\sf S$, $\sf P$, $\sf P_U$, $\sf E_{el}$ be the homomorphic image, subalgebra, direct product, ultraproduct and elementary extension class operators. 

\begin{proposition}[\protect{\cite[Theorem V.2.25]{BS81}\cite[Corollary 2.3.4]{Gor98}}]\label{prop::quasivariety generation}
For a class $\mathcal K$ of algebras in the same language we have ${\sf Q}(\mathcal K)={\sf SPP_U}(\mathcal K)$. In particular, for an algebra $\A$ we have ${\sf Q}(\A)={\sf SPE_{el}}(\A)$.
\end{proposition}

Let $\Q$ be a quasivariety. A congruence
$\alpha$ on an algebra $\A$ is called a \mbox{\emph{$\Q$-congruence}}
provided $\A/\alpha\in\Q$. Note that $\A\in\Q$ if and only if the equality relation on $A$ is a $\Q$-congruence of $\A$. The set $\Con_\Q(\A)$ of all
\mbox{$\Q$-congruences} of $\A$ forms an algebraic lattice. It 
is a meet-subsemilattice of $\Con(\A)$, the lattice of all congruences of $\A$ \mbox{\cite[Corollary 1.4.11]{Gor98}}.

A nontrivial algebra $\S$ is $\Q$-simple if $\Con_\Q(\S)$ has exactly two elements: the equality relation $\id_S$ on $S$ and the total relation $S^2$ on $S$. A nontrivial algebra $\S\in\Q$ is
$\Q$-\emph{subdirectly irreducible} if the equality relation on $A$ is completely meet irreducible in $\Con_\Q(\A)$. (In the case when $\Q$ is a variety we drop the prefix \mbox{``$\Q$-''.)} The importance of $\Q$-subdirectly irreducible algebras follows from the fact that they determine $\Q$. Indeed, in an algebraic lattice each element
is a meet of completely meet-irreducible elements. Thus we have the following fact.

\begin{proposition}[\protect{\cite[Theorem 3.1.1]{Gor98}}]\label{prop::Mal'cev thm}
Every algebra in $\Q$ is isomorphic to a subdirect product of $\Q$-subdirectly irreducible algebras. In particular, the class of $\Q$-subdirectly irreducible algebras generates $\Q$.
\end{proposition}

A quasivariety $\Q$ is \emph{semisimple} if every $\Q$-subdirectly irreducible algebra is $\Q$-simple. So if $\Q$ is semisimle, then it is generated by $\Q$-simple algebras.

Let $\G\in\Q$ and $W\subseteq G$. We say that $\G$ is \emph{free for $\Q$ over $W$}, and is \emph{of rank} $|W|$, if $\G\in\Q$ and it satisfies the following \emph{universal mapping property}: every mapping $f\colon W\to A$, where $A$ is a carrier of an algebra $\A$ in $\Q$, is uniquely extendable to a homomorphism $\bar f\colon\G\to\A$. In our consideration the following weaker property will be used. If $\G$ is free for $\Q$ and $\A\in\Q$, then there is a homomorphism from $\G$ into $\A$ (we adopt the convention that algebras have nonempty carriers). 
If $\Q$ contains a nontrivial algebra, then it has free algebras over arbitrary non-empty sets and, in fact, they coincide with free algebras for the variety $\sf V(\Q)$. (Note here that $\sf V(\Q)$ is the class of all homomorphic images of algebras from $\Q$.) Let us fix a denumerable set of variables $V=\{v_0,v_1,v_2,\ldots\}$. We denote a free algebra for $\Q$ over $V$ by $\Fr$ and the free algebra for $\Q$ over $V_k=\{v_0,v_1,\ldots,v_{k-1}\}$ by $\Fr(k)$. One may construct $\Fr$ and $\Fr(k)$ by forming the algebra of terms over $V$, or $V_k$ respectively, and divide it by a certain congruence. This congruence identifies terms $s(\bar v),t(\bar v)$ which determine the same term operation on every algebra from $\Q$ (in other words, when $\Q\models(\forall\bar x)[t(\bar x)\=s(\bar x)]$). We will notationally identify
terms with the elements of $\Fr$ that they represent.

For an algebra $\A$ and a set $H\subseteq A^2$ there exists the least $\Q$-con\-gruence $\theta_\Q(H)$ on $\A$ containing $H$.
We say that an algebra is \emph{$\Q$-finitely presented}
if it is isomorphic to an algebra of the form $\Fr(k)/\theta_\Q(H)$ for some natural number $k$ and some finite set $H$ \cite[Chapter 2]{Gor98}. For a tuple $\bar x=(x_0,\ldots,x_{k-1})$ of variables and a conjunction of equations $\varphi (\bar x)=s_1(\bar x)\=t_1(\bar x)\meet\cdots\meet s_n(\bar x)\=t_n(\bar x)$ let
\[
\P_{\varphi(\bar x)}=\Fr(k)/\theta_\Q(\{(s_1(\bar v),t_1(\bar v)),\ldots,(s_n(\bar v),t_n(\bar v))\}),
\]
where $\bar v=(v_0,\ldots,v_{k-1})$. Note that every finitely presented algebra is isomorphic to some $\P_{\varphi(\bar x)}$.
In the following fact we notationally identify variables from $\bar v$ with their congruence classes.
\begin{fact}\label{fact::fp-algebras}
\mbox{}
\begin{enumerate}
\item A quasivariety $\Q$ satisfies a quasi-identity $(\forall \bar x)[\varphi(\bar x)\to\psi(\bar x)]$ if and only if $\P_{\varphi(\bar x)}\models \psi(\bar v)$.
\item For a conjunction of equations $\varphi(\bar x)$, for every algebra $\A\in\Q$ there exists a homomorphism from $\P_{\varphi(\bar x)}$ into $\A$ if and only if $\A\models(\exists\bar x)\, \varphi(\bar x)$.
\end{enumerate}
\end{fact}

\section{General solution}

Let $\Q$ be a quasivariety and $\Fr$ be a free algebra of denumerable rank for $\Q$. A quasi-identity which is true in $\Fr$ is called \emph{$\Q$-admissible}. For a quasi-identity $q=(\forall \bar x)\, [\varphi(\bar x)\to \psi(\bar x)]$ let
\[
q^*=(\forall \bar x)\, [\neg\varphi(\bar x)].
\]
We say that a $\Q$-admissible quasi-identity $q$ is \emph{$\Q$-passive} if $\Fr\models q^*$, and \emph{$\Q$-active} otherwise. By Fact \ref{fact::fp-algebras}, a $\Q$-admissible quasi-identity $q$ is $\Q$-active iff $\P_{\varphi(\bar x)}$ admits a homomorphism into $\Fr$, where $\varphi(\bar x)$ is the premise of $q$.

A quasivariety $\Q$ is \emph{structurally complete} (SC) provided that every $\Q$-admissible quasi-identity is true in $\Q$. In other words, if $\Q={\sf Q}(\Fr)$.
A quasivariety $\Q$ is \emph{almost structurally complete} (ASC) provided that every $\Q$-active quasi-identity holds in $\Q$.

\begin{lemma}\label{lem:: homo into el ext of F}
Let $\Q$ be an SC quasivariety. Then
every nontrivial algebra $\B$ from $\Q$ admits a homomorphism into some elementary extension $\G$ of $\Fr$. 
\end{lemma}

\begin{proof}
By Proposition \ref{prop::quasivariety generation} and the assumption that $\Q$ is SC, the algebra $\B$ belongs to ${\sf SPE_{el}}(\Fr)$. This means that for each pair of  two  distinct elements $a$ and $b$ from $B$ there is a homomorphism $h\colon\B\to\G$, where $\G\in{\sf E_{el}}(\Fr)$, separating $a$ and $b$. Since $\B$ is nontrivial, there exists at least one such homomorphism.
\end{proof}

\begin{theorem}\label{thm::general solution}
Let $\Q$ be an ASC quasivariety. Then the following conditions are equivalent
\begin{enumerate}
\item $\Q$ is SC;
\item every nontrivial algebra from $\Q$ admits a homomorphism into some elementary extension of $\Fr$;
\item every nontrivial $\Q$-finitely presented algebra admits a homomorphism into $\Fr$.
    \end{enumerate}
\end{theorem}

\begin{proof}\mbox{}\\
\noindent(1)$\Rightarrow$(2)  It follows from Lemma \ref{lem:: homo into el ext of F}.

\noindent(2)$\Rightarrow$(3) It follows from Fact \ref{fact::fp-algebras} Point (2).

\noindent(3)$\Rightarrow$(1) Let $q=(\forall\bar x)[\varphi(\bar x)\to\psi(\bar x)]$ be a $\Q$-admissible quasi-identity. We want to show that $\Q\models q$. Since $\Q$ is ASC, we may assume that $q$ is $\Q$-passive. This exactly means that $\P_{\varphi(\bar x)}$ does not admit a homomorphism into $\Fr$. Then Condition (3) says that $\P_{\varphi(\bar x)}$ is trivial. Hence, by Fact \ref{fact::fp-algebras} point (1), $\Q\models q$ (whatever $\psi$ is).
\end{proof}

\begin{corollary}\label{cor:: idempotent in G}
If $\Q$ is an ASC quasivariety such that $\Fr$ has an elementary extension $\G$ with an idempotent element, then $\Q$ is SC.
\end{corollary}

\section{Semisimple quasivarieties}\label{sec::semisimple}

Let us now prove a lemma which is built upon \cite[Corollary 2.8]{Ber88}.

\begin{lemma}\label{lem::key lemma}
Let $\Q$ be an SC quasivariety such that every elementary extension of $\Fr$ does not have idempotent elements. If $\S$ is a $\Q$-simple algebra and $\A$ is a nontrivial algebra from $\Q$, then $\S$ embeds into some elementary extension of $\A$.
\end{lemma}

\begin{proof}
By Lemma \ref{lem:: homo into el ext of F}, the lack of idempotent elements in every elementary extension of $\Fr$ yields the lack of idempotent elements in every nontrivial algebra $\B\in\Q$. Thus, by the $\Q$-simplicity of $\S$ every homomorphism from $\S$ into such $\B$ is injective. It follows that we only need to verify the existence of any homomorphism from $\S$ into some elementary extension of $\A$.
This will be accomplished by showing that the set $\diagplus(\S)\cup\diagel(\A)$ is satisfiable.
Here $\diagplus(\S)$ is the positive diagram of $\S$, i.e., the set of atomic sentences which are valid in the expansion of $\S$ obtained by adding constants corresponding to all elements of $S$. And $\diagel(\A)$ is the elementary diagram of $\A$, i.e., the set of first order sentences which are valid in the analogous expansion of $\A$ \cite{CK90} \cite{Hod93}.

By the compactness theorem, it is enough to show that for every finite $\Sigma\subseteq\diagplus(\S)$ the set $\Sigma\cup\diagel(\A)$ is satisfiable. Let $\sigma(\bar x)$ be the formula obtained from the conjunction $\bigmeet\Sigma$ by replacing every constant which is not a constant of $\S$ with a new variable. Then the satisfiability of $\Sigma\cup\diagel(\A)$ is equivalent to the satisfiability of $\{\sigma(\bar x)\}\cup\diagel(\A)$. The latter will be shown by the verification that $\A\models(\exists \bar x)\,\sigma(\bar x)$.

By Lemma \ref{lem:: homo into el ext of F}, there is a homomorphism from $\S$ into some elementary extension $\G$ of $\Fr$. Since $\sigma(\bar x)$ is positive and $\S\models(\exists \bar x)\,\sigma(\bar x)$, we have $\G\models(\exists \bar x)\,\sigma(\bar x)$ and also $\Fr\models(\exists \bar x)\,\sigma(\bar x)$. By the freeness of $\Fr$ there also exists a homomorphism from $\Fr$ into $\A$. Thus, using once more the positivity of $\sigma(\bar x)$ we obtain $\A\models(\exists \bar x)\,\sigma(\bar x)$.
\end{proof}

At this point we are able to characterize SC quasivarieties among semisimple ASC quasivarieties.

\begin{theorem}\label{thm:: main}
Let $\Q$ be an ASC semisimple quasivariety. Then $\Q$ is SC if and only if there exists an elementary extension of $\Fr$ with an idempotent element or $\Q$ is minimal.
\end{theorem}

\begin{proof}
If there exists an elementary extension of $\Fr$ with an idempotent element, Corollary \ref{cor:: idempotent in G} implies that $\Q$ is SC. If $\Q$ is minimal, then it has exactly two subquasivarieties: the trivial one and itself. One of them is $\sf Q(\Fr)$ and it cannot be the trivial one since $\Fr$ is nontrivial. Hence then $\Q$ is also SC.

For the converse implication assume that $\Q$ is an SC quasivariety such that every elementary extension of $\Fr$ does not have any idempotent element. By Lemma \ref{lem::key lemma}, every $\Q$-simple algebra embeds into some elementary extension of $\Fr$. In particular, all of them belong to $\sf Q(\Fr)$. Thus, by the semisimplicity, Proposition \ref{prop::Mal'cev thm} yields that $\Q=\sf Q(\Fr)$ and $\Q$ is SC.
\end{proof}

When an additional finiteness condition is imposed the formulation of Theorem \ref{thm:: main} may be simplified. Note however that Example \ref{exm:: idemp w G} shows that such simplification is not possible in general.

\begin{corollary}\label{cor:: SC with finiteness condition}
Let $\Q$ be an ASC semisimple quasivariety. Assume that either the language of $\Q$ is finite or $\Fr$ has a subalgebra $\A$ with the finite carrier $A$. Then $\Q$ is SC if and only if $\Fr$ has an idempotent element or $\Q$ is minimal.
\end{corollary}

\begin{proof}
By Theorem \ref{thm:: main}, it is enough to show that if there exists an elementary extension $\G$ of $\Fr$ with an idempotent element, then $\Fr$ also has an idempotent element.

Assume first that the language of $\Q$ is finite. Then the existence of an idempotent element is expressible by a first order sentence. Since $\G$ and $\Fr$ satisfy the same first order sentences,  $\Fr$ has an idempotent element.

Now assume that $\A$ is a subalgebra of $\Fr$ with the finite carrier. By \cite[Theorem V.2.16]{BS81}, $\G$ elementarily embeds into some ultrapower $\Fr^I/U$. Thus there is an idempotent element $e$ in $\Fr^I/U$. Let $f\colon\Fr\to\A$ be any homomorphism. Its existence is guaranteed by the freeness of $\Fr$. We have also a homomorphism $f^*\colon\Fr^I/U\to\A^I/U$ given by $t/U\mapsto f(t)/U$. Then $f^*(e)$ is idempotent in $\A^I/U$. But the finiteness of $\A$ yields that $\A$ and $\A^I/U$ are isomorphic \cite[Theorem IV.6.5]{BS81}. Hence $\A$ and $\Fr$ have an idempotent element.
\end{proof}

\section{Discriminator varieties}

The ternary operation on a set given by 
\[
\disc(a,b,c)=
\begin{cases}
a & \text{ if } a\neq b\\
c & \text{ if } a= b
\end{cases}
\]
is called the \emph{discriminator operation}.
A variety $\V$ is a \emph{discriminator variety} if it is generated by a class $\mathcal K$ of algebras for which there is a term whose interpretation is the discriminator operation in every member of $\mathcal K$.
The following facts reveals why the results from Section \ref{sec::semisimple} are applicable for discriminator varieties.

\begin{proposition}\label{prop::discriminator properties}
Let $\V$ be a discriminator variety generated by a class $\mathcal K$.  Assume that there is a term whose interpretation is the discriminator operation in every member of $\mathcal K$.  Then
\begin{enumerate}
\item $\V$ is ASC;
\item $\V$ is semisimple;
\item A nontrivial algebra from $\V$ is simple if and only if it belongs to ${\sf SP_U}(\mathcal K)$;
\end{enumerate}
\end{proposition}

\begin{proof}
For (1) note that $\V$ has projective unification \cite[Theorem 3.1]{Bur92} and this yields ASC \cite[Corollary 5.1]{DS14}.
For (2) and (3) see \cite[Theorem IV.9.4]{BS81} or \cite[Lemma 1.3]{BW77}.
\end{proof}

A relevant fact about discriminator varieties is that if they are minimal (as varieties) they are also minimal as quasivarieties. As far as we know the general case is a new result. However,
this was already proved in \cite{BCV01} under some fairly general additional hypothesis. It includes the cases when the language is finite or the variety is locally finite. Note also that in the locally finite case even congruence-modularity yields that minimal varieties are minimal as quasivarieties \cite[Corollary 13]{BM90}. In \cite[Examples 14 and 15]{BM90} the reader may also find examples of minimal varieties which are not minimal as quasivarieties. One of them is locally finite and another is an arithmetical variety.    

In the proof we make use of Boolean products, and some of their properties. Let us recall here that Boolean products are some special subdirect products. For the precise definition of Boolean product and the proof of the fundamental theorem below we refer the reader to e.g. \cite[Theorem 9.4]{BS81}. 

\begin{theorem}\label{thm::bu-fl ke we}
Every algebra in a discriminator variety is a isomorphic with a Boolean product of trivial and simple factors.
\end{theorem}


We shall also need the following preservation result.

\begin{lemma}[\protect{\cite[Lemma 9.3]{BW79}}]\label{thm::pp preservation}
Boolean products preserve the satisfaction of primitive positive sentences, i.e., sentences of the form $(\exists \bar x)\,[\bigwedge p_{i}(\bar x))\=q_{i}(\bar x)]$. This means that if such a sentence is true of every factor then it is true of the Boolean product.
\end{lemma}

\begin{proposition}\label{thm::minimal discr}
If $\mathcal{V}$ is a minimal discriminator variety, then $\mathcal{V}$ is a
minimal quasivariety.
\end{proposition}

\begin{proof}
Let $q$ be a quasi-identity and suppose that there is a nontrivial algebra $\A$ from $\V$ satisfying $q$. Our aim is to show that $\V\models q$. By \cite[Theorem 2.4]{Vag95}, 
we have that there are terms $s,t,u,v$ such that
\[
\V\models q \leftrightarrow \left [(\forall \bar x)\, s(\bar x)\=t(
\bar x)\vee (\forall \bar x)\, u(\bar x)\not \= v(\bar x)\right].
\]
Since $\mathcal{V}$ is a minimal variety, either
\[
\mathcal{V}\models (\forall \bar x)\, s(\bar x)\=t(\bar x),
\]
and then clearly $\V\models q$, or the subvariety of $\V$ defined by ($\forall \bar x)\, s(\bar x)\=t(\bar x)$ is trivial. Thus in the latter case we  
have 
\[
\A\models (\forall \bar x)\, u(\bar x)\not \= v(\bar x).
\]
In view of Theorem \ref{thm::bu-fl ke we}, we can suppose that $\A\leq \prod \{\A_{i}:i\in I\}$ is
a Boolean product, where each $\A_{i}$ is simple or trivial. Now Lemma \ref{thm::pp preservation} yields that there exists $i\in I$ such that $$\A_{i}\models (\forall \bar x)\, u(\bar x)\not \= u(\bar x).$$ 
Note that $\A_i$ is non-trivial and hence $\A_i$ is simple.
Since $\mathcal{V}$ is a minimal variety it is generated by $\A_i$. Thus by (3) of Proposition \ref{prop::discriminator properties}, we have
\[
\V_{S}\subseteq {\sf SP_U}(A_i),
\]
where $\V_{S}$ denotes the class of simple algebras in $\mathcal{V}$. Hence
\[
\V_{S}\models (\forall \bar x)\, u(\bar x)\not \= v(\bar x)
\]
since the satisfaction of $(\forall \bar x)\, r(\bar x)\not \= s(\bar x)$ is preserved by 
 $\sf S$ and $\sf P_U$ class operators. Since it is also preserved by forming nontrivial subdirect products, Proposition \ref{prop::Mal'cev thm} implies that
\[
\V - \{\text{trivial algebras}\}\models (\forall \bar x)\, u(\bar x)\not \= v(\bar x).
\]
This finally gives that $\V\models q$.
\end{proof}

\begin{theorem}\label{thm:: disr var}
Let $\V$ be a discriminator variety. Then $\V$ is SC if and only if there exists an elementary extension of $\Fr$ with an idempotent element or $\V$ is a minimal (as a variety or quasivariety).
\end{theorem}

\begin{proof}
It follows from Theorem \ref{thm:: main}, points (1) and (2) of Proposition 
\ref{prop::discriminator properties} and Proposition \ref{thm::minimal discr}.
\end{proof}

The following corollary shows the advantage of dealing with ASC comparied to SC. 
It is a quite common situation in algebraic logic when there is no idempotent element in any extension of $\Fr$.

\begin{corollary}\label{cor:: disc with constants}
Let $\V$ be a discriminator variety. Assume that there are two distinct constants in $\Fr$. Then $\V$ is SC if and only if it is minimal.
\end{corollary}

\begin{proof}
Indeed, an idempotent element, by the definition, is a one-element subalgebra. But the smallest subalgebra of $\Fr$ has at least two elements. 
\end{proof}

\begin{example}
Let us consider the variety of two-dimensional (representable) cylindric algebras $\mathcal{(R)CA}_2$ or the variety two-dimensional representable cylindric algebras $\mathcal{RCA}_2$ \cite{HMT85}. It is known that it is a discriminator variety with continuum many subvarieties \cite[Theorem 4.2]{Bez04}. But is has only two minimal subvarieties \cite[Corollary 7.5]{Bez04}. Thus, by Corollary \ref{cor:: disc with constants} and Proposition \ref{prop::discriminator properties} Point (1), there are continuum many subvarieties of $\V$ which are ASC but are not SC.
\end{example}

\begin{example}
We have similar situation for relation algebras.
(Representable) relation algebras form a discriminator variety. It has exactly three minimal subvarieties \cite{Tar56}. There exists continuum many varieties of (representable) relation algebras \cite[Theorem 7.8]{Jon82}. All of them are ASC.
\end{example}

By imposing the same finiteness condition as in Corollary \ref{cor:: SC with finiteness condition} we again obtain a certain simplification.

\begin{corollary}\label{cor:: discriminator+finite condition}
Let $\V$ be a discriminator variety. Assume that either the language of $\V$ is finite or $\Fr$ has a finite subalgebra. Then $\V$ is SC if and only if $\Fr$ has an idempotent element or $\V$ is a minimal variety.
\end{corollary}


\begin{example}\label{exm:: idemp w G}
We will construct an SC discriminator variety $\V$ which is not minimal (as a variety and hence also as a quasivariety) and its free algebras have no idempotent element.

Let $\W$ be the variety in the language with denumerable many unary operations $f_i$, $i\in\Nat$, defined by the identities
\begin{gather*}
(\forall x)\,f_i(f_j(x))\=f_j(f_i(x))\\
(\forall x)\,f_i(f_i(x))\=f_i(x)
\end{gather*}
for all $i,j\in\Nat$. Let $\A$ be the algebra in $\W$ with the carrier $A=\{a,b\}$ and such that $f_i(a)=b$ for all $i\in\Nat$. Let $\B$ be a free algebra for $\W$ of rank at least 1. Let us expand the language of $\W$ by adding one ternary function symbol $\disc$. Let $\A_{\disc}$ and $\B_{\disc}$ be the expansions of $\A$ and $\B$ respectively to the new language where $\disc$ is interpreted as the discriminator operation. Let $\V$ be the variety generated by $\A_{\disc}$ and $\B_{\disc}$. Then $\V$ is a discriminator variety. Each free algebra for $\V$ does not have any idempotent element since $\B_{\disc}$ does not. Also $\V$ is not minimal. In order to see this note that $\B$ is simple and does not belong to $\sf S(\A)$. Since $\A$ is finite, $\sf SP_U(\A)= S(\A)$. Then Proposition \ref{prop::discriminator properties} Point (3) yields that $\B_{\disc}\not\in {\sf V}(\A_{\disc})$. 

But $\Fr$ has an ultrapower with an idempotent element, and hence is SC. To see this let us take an element $u\in F$ and consider a sequence $\bar u$ in $F^\Nat$ given by $i\mapsto f_{i-1}(f_{i-2}(\cdots f_0(u)\cdots))$. Now if $U$ is a non-principal
ultrafilter over $\Nat$, then $\bar u/U$ is idempotent in $\Fr^\Nat/U$.
\end{example}

\begin{example} \label{exm:: minimal discriminator without homo into F}
Here we will present a minimal discriminator $\V$ variety with a countably algebra $\S$ which does not admit a homomorphism into $\Fr$. Recall that by Theorem \ref{thm:: disr var}, $\V$ is SC. Thus, by Theorem \ref{thm::general solution}, every $\V$-finitely presented algebra admits a homomorphism into $\Fr$.  

Let $\A=(\Nat,s,\disc)$ be the algebra, with the set of natural numbers as its carrier and with the successor  and the discriminator operations on $\mathbb N$ as the basic operations. Let $\V=\sf V(\A)$. By Proposition \ref{prop::discriminator properties} point (3), every simple algebra in $\V$ has a subalgebra isomorphic to $\A$. Thus $\V$ is minimal and SC. For $\S$ we may take a simple algebra $(\mathbb Z,s,\disc)$ defined similarly as $\A$ but with the set of integers as the carrier. Note that $\S$  is isomorphic to the subalgebra $\A^\Nat/U$, where $U$ is a non-principal ultrafiler on $\Nat$, generated by the identity function on $\Nat$. 

And indeed, the algebra $\S$ does not admit a homomorphism into any free algebra for $\V$. This fact follows from the simplicity of $\S$, the lack of an idempotent element in every free algebra $\G$ for $\V$, and the following facts.

\begin{claim}
The operation $s$ is injective in $\G$.
For every $a\in G$, there exist $n\in\mathbb N$ and $b\in G$ such that $a=s^n(b)$ and $b\not\in s(G)$.
\end{claim}

\begin{proof}
Free algebras for $\V$ are isomorphic to subalgebras of powers of $\A$. For instance, $\Fr(1)$ is isomorphic to the subalgebra of $\A^A$ generated by the identity function on $A$. So it is enough to verify the property from the claim for $\G$ which is a subalgebra of $\A^I$ for some $I$. 

The injectivity of $s$ in $\G$ follows from the injectivity of $s$ in $\A$. Let $n$ be the smallest number $k$ such that the function $i\mapsto a(i)-k$, $i\in I$, belongs to $G$, and take $b\colon i\mapsto a(i)-n$ for $i\in I$.
\end{proof}

\end{example}



\section{Conclusion}

We solved the SC vs ASC problem in the case of semisimple quasivarieties. In particular, we obtained a characterization of SC discriminator varieties. Our results shows that SC property may be very restrictive compared to ASC. It is a common situation that for a discriminator variety the only SC subvarieties are minimal ones while all of them are ASC. In our opinion, it shows the advantage of dealing with ASC instead of SC.

\bibliographystyle{plain}
\bibliography{scdisc03082014}

\end{document}